\documentclass[graybox]{article}

\usepackage{type1cm}        
\usepackage{graphicx}        
\usepackage{multicol}        
\usepackage[bottom]{footmisc}

\usepackage{newtxtext}       %
\usepackage{newtxmath}       
 \usepackage[backref,colorlinks,linkcolor=red,anchorcolor=green,citecolor=blue]{hyperref}

 \newtheorem{theorem}{Theorem}

\begin{document}

\title{A review on a general multi-species BGK model: modelling, theory and numerics}
\author{Marlies Pirner and Sandra Warnecke}

%
%
\maketitle
\footnote{Marlies Pirner \\ Julius-Maximilians-Universität Würzburg, Emil-Fischer-Straße 40, 97074 Würzburg, \\ marlies.pirner@mathematik.uni-wuerzburg.de \\
 Sandra Warnecke \\ Julius-Maximilians-Universität Würzburg, Emil-Fischer-Straße 40, 97074 Würzburg, \\ sandra.warnecke@mathematik.uni-wuerzburg.de}
\abstract{In this article we focus on kinetic equations for gas mixtures since in applications one often has to deal with mixtures instead of a single gas. In particular we consider an approximation of the Boltzmann equation, the Bathnagar-Gross-Krook (BGK) equation. This equation is used in many applications because it is very efficient in numerical simulations. In this article, we recall 
a general BGK equation for gas mixtures which has free parameters. Specific choices of these free parameters lead to special cases in the literature.  For this model, we provide an overview concerning modelling, theoretical results and numerics.}

\section{Introduction}
 In this paper we shall concern ourselves with a kinetic description of gas mixtures. For simplicity in notation and statements, we present it here for two species, but the model can be extended to  an arbitrary number of species since we only consider binary interactions.  A gas of mono atomic molecules and two species is traditionally described via the Boltzmann equation for the distribution functions $f_1=f_1(x,v,t), f_2=f_2(x,v,t)$, see for example \cite{Cercignani, Cercignani_1975}. Here, $x \in \mathbb{R}^3$ and $v \in \mathbb{R}^3$ are the phase space variables, position and velocity of the particles, and $t \geq 0$ denotes the time. \textcolor{black}{Assume that the particles of species $1$ have mass $m_1$ and the particles of species $2$ have mass $m_2$.} The Boltzmann equation for gas mixtures is of the form
 \begin{align*}
     \partial_t f_1 + v \cdot \nabla_x f_1 = Q_{11}(f_1,f_1) + Q_{12}(f_1,f_2),\\ \partial_t f_2 + v \cdot \nabla_x f_2 = Q_{22}(f_2,f_2) + Q_{21}(f_2,f_1), 
 \end{align*}
 where the intra-species collision operators $Q_{11}(f_1,f_1)$ and$Q_{22}(f_2,f_2)$ satisfy 
 \begin{align} 
\int Q_{kk}(f_k,f_k) \begin{pmatrix} 1 \\ m_k v \\ m_k |v|^2 \end{pmatrix} dv = 0, \quad \quad k=1,2 
\label{cons1}
\end{align}
and the inter-species collision operators  $Q_{12}(f_1,f_2)$ and $Q_{21}(f_2,f_1)$ satisfy
\begin{align}
\begin{split}\label{cons2}
\int Q_{12}(f_1,f_2) dv = \int Q_{21}(f_2,f_1) dv = 0, \\ \int \left( \begin{pmatrix} m_1 v \\ m_1 |v|^2 \end{pmatrix}  Q_{12}(f_1,f_2)+ \begin{pmatrix} m_2 v \\ m_2 |v|^2 \end{pmatrix} Q_{21}(f_2,f_1)\right) dv = 0. 
\end{split}
\end{align}
These properties of the collision operator ensure conservation of the number of particles, total momentum and total energy at the macroscopic level, \eqref{cons1} in intra-species interactions, \eqref{cons2} in inter-species interactions.
In addition, the collision operators satisfy the inequalities
\begin{align}
\begin{split}
    \int Q_{kk}(f_k,f_k) \ln f_k dv \leq 0, \quad \quad k=1,2
\\
    \int Q_{12}(f_1,f_2) \ln f_1 dv + \int Q_{21}(f_2,f_1) \ln f_2 dv \leq 0.
    \end{split}
    \label{entropy}
\end{align}
The first inequality turns into an
equality if and only if $f_k$ is a Maxwell distribution $M_k$ given by 
\begin{align} 
M_k &= \frac{n_k}{\textcolor{black}{(2 \pi \frac{T_k}{m_k})^{3/2}} }  \exp\left({- \frac{|v-u_k|^2}{2 \frac{T_k}{m_k}}}\right).
\label{Max_one}
\end{align}
 Here we define for any $f_1,f_2: \Lambda \subset \mathbb{R}^3 \times \mathbb{R}^3 \times \mathbb{R}^+_0 \rightarrow \mathbb{R}$ with $(1+|v|^2)f_1,(1+|v|^2)f_2 \in L^1(\mathbb{R}^3), f_1,f_2 \geq 0$\textcolor{black}{,} the  macroscopic quantities 
\begin{align}
\int f_k(v) \begin{pmatrix}
1 \\ v  \\ m_k |v-u_k|^2 
\end{pmatrix} 
dv =: \begin{pmatrix}
n_k \\ n_k u_k \\ 3 n_k T_k
\end{pmatrix} , \quad k=1,2
\label{moments}
\end{align} 
where $n_k$ is the number density, $u_k$ the mean velocity and $T_k$ the mean temperature of species \textcolor{black}{ $k$ ($k=1,2$).} For ease we write $T_k$ instead of $k_B T_k$, where $k_B$ is Boltzmann's constant. 

In the second inequality in \eqref{entropy}, we have equality if and only if $f_1$ and $f_2$ are Maxwell distributions $M_1$ and $M_2$ and additionally if and only if $u_1=u_2$ and $T_1=T_2$.

 \textcolor{black}{If we are close to equilibrium,} the complicated interaction terms of the Boltzmann equation can be simplified by a so called BGK approximation, consisting of a collision frequency $\nu_{kj} n_j$ multiplied by the deviation of the distributions from a local Maxwell distribution  
 \begin{align} \begin{split} \label{BGK}
\partial_t f_1 + v \cdot \nabla_x  f_1   &= \nu_{11} n_1 (M_1-f_1) + \nu_{12} n_2 (M_{12}- f_1),
\\ 
\partial_t f_2 + v \cdot \nabla_x  f_2 &=\nu_{22} n_2 (M_2 - f_2) + \nu_{21} n_1 (M_{21}- f_2).
\end{split}
\end{align}
The collision frequencies per density $\nu_{kj}$ are assumed to be dependent only on $x$ and $t$ and not on the microscopic velocity $v$. For references taking into account also a dependency on the microscopic velocity $v$ see \cite{Struchtrupp} for the one species case, \cite{Pirner_velocity} for the gas mixture case and \cite{Sandra_Numerik} for the numerics of the gas mixture case.

The mixture Maxwell distributions $M_{12}$ and $M_{21}$ are given by
\begin{align}
\begin{split}
M_{12}(x,v,t) = \frac{n_{12}}{(2 \pi \frac{T_{12}}{m_1})^{3/2} }  \exp\left({- \frac{|v-u_{12}|^2}{2 \frac{T_{12}}{m_1}}}\right), \\
M_{21}(x,v,t) = \frac{n_{21}}{(2 \pi \frac{T_{21}}{m_2})^{3/2} }  \exp\left({- \frac{|v-u_{21}|^2}{2 \frac{T_{21}}{m_2}}}\right).
\end{split}
\label{BGKmix}
\end{align}
where $n_{kj}, u_{kj}$ and $T_{kj}$ will shortly be defined.

 This approximation should be constructed in a way such that it  has the same main properties as the Boltzmann equation mentioned above. 


Now, the question arises how to choose the mixture quantities $n_{kj}, u_{kj}$ and $T_{kj}$. In this review article we present a general model which is published in \cite{Pirner} for two species. This model  contains a lot of proposed models in the literature as special cases.
Examples are the models of Gross and Krook \cite{gross_krook1956}, Hamel \cite{hamel1965},  Asinari \cite{asinari}, Garzo Santos Brey \cite{Garzo1989}, Sofena \cite{Sofonea2001}, Cercignani \cite{Cercignani_1975}, Greene \cite{Greene} and recent models by Bobylev, Bisi, Groppi, Spiga, Potapenko \cite{Bobylev}; Haack, Hauck, Murillo \cite{haack}. 

The second last \cite{Bobylev} presents an additional motivation how the corresponding model can be derived formally from the Boltzmann equation, whereas the last \cite{haack} presents a derivation to macroscopic equations on the Navier-Stokes level and numerical results. 

 BGK  models give rise to efficient numerical computations, which are asymptotic preserving, that is they remain efficient even approaching the hydrodynamic regime \cite{Puppo_2007, Jin_2010,Dimarco_2014, Bennoune_2008,  Bernard_2015, Crestetto_2012}. \textcolor{black}{ However, the BGK approximation is incapable of reproducing the correct Boltzmann hydrodynamic regime in the asymptotic continuum limit. Therefore, a modified version called the ES-BGK approximation was suggested  by Holway for one species \cite{Holway}. Then the H-Theorem of this model  was shown in \cite{Perthame} and existence and uniqueness of  mild solutions in \cite{Yun_mild}. Alternatively, the Shakov model \cite{Shakov} and a BGK model with velocity dependent collision frequency \cite{Struchtrupp} were suggested to achieve the correct Prandtl number. For the BGK model with velocity dependent collision frequency, it is shown that a power law for the collision frequency also leads to the proper Prandtl number.  
 The standard BGK model is extended to a velocity dependent collision frequency while it still satisfies the conservation properties. This works when replacing
 the Maxwell distribution by a different function, for details see \cite{Struchtrupp}. For this model an H-Theorem can be proven. The existence of these modified functions is proven in \cite{Pirner_velocity}. However, since BGK models form the basis to build extended models as ES-BGK models, Shakov models and BGK models with velocity dependent collision frequency, we will mainly review BGK models for gas mixtures of the form \eqref{BGK} in this paper.}

\textcolor{black}{Considerations of the hydrodynamic regime for an BGK model for gas mixtures is considered for example in \cite{haack}, a special case of the model presented in this paper. It presents a Chapman-Enskog expansion  with transport coefficients in section 5, a comparison with other BGK models for gas mixtures in section 6  and a numerical implementation.} 

Additionally, we want to mention that there is also another type of BGK model for gas mixtures containing only  one collision term on the right-hand side. Examples for this are Andries, Aoki and Perthame \cite{AndriesAokiPerthame2002} and the models in \cite{Brull_2012, Groppi}. \textcolor{black}{A derivation of the Navier-Stokes system in the compressible regime for the model in \cite{AndriesAokiPerthame2002} and the corresponding transport coefficients can be found in section 4 of \cite{AndriesAokiPerthame2002}.} \textcolor{black}{The transport coefficients of the hydrodynamic regime for the model in \cite{Brull_2012} can be found in section 5 of} \cite{Brull_2012}. \textcolor{black}{A comparison of these models concerning their hydrodynamic limit can be found in \cite{Boscarino}.}
\textcolor{black}{For gas mixtures there are also many results concerning  extensions to ES-BGK models, Shakov models and BGK models with velocity dependent collision frequency \cite{Groppi,Brull,Blaga2,Pirner_velocity}.}

In the following, we will present theoretical and numerical results for this  general BGK model for two species with two interaction terms which captures all those special cases in the literature.
The outline of the paper is as follows: In section \ref{sec2} we will present the general multi-species BGK model for two species. For this model, we will give  a review of recent theoretical results in section \ref{sec3}. The physical meaning and possible choices of the free parameters are discussed in section \ref{sec:parameters}.  And recent existing numerical schemes are given in section \ref{sec4}. 

\section{The general BGK model for gas mixtures}
\label{sec2}

In this section, we will concern the question of how to choose the mixture quantities $n_{12}, n_{21}, u_{12}, u_{21}, T_{12}, T_{21}$ and the collision frequencies. The collision frequencies $\nu_{11} n_1$ and $\nu_{22} n_2$ correspond to interactions 
of the particles of each species with itself, while $\nu_{12} n_2$ and $\nu_{21} n_1$ are related to inter-species collisions.
To be flexible in choosing the relationship between the collision frequencies, we now assume 
\begin{equation} 
\nu_{12}=\varepsilon \nu_{21}, \quad 0 < \varepsilon \leq 1.
\label{coll}
\end{equation}
The restriction on $\varepsilon$ is without loss of generality. If $\varepsilon >1$, exchange the notation $1$ and $2$ and choose $\frac{1}{\varepsilon}.$ In addition, we assume that all collision frequencies are positive. 
The Maxwell distributions $M_1$ and $M_2$ in \eqref{Max_one} are chosen to have the same density, mean velocity and temperature as $f_1$ and $f_2$, respectively. With this choice, we guarantee the conservation of mass, momentum and energy in interactions of one species with itself \eqref{cons1} (see section 2.2 in \cite{Pirner}).
The remaining parameters $n_{12}, n_{21}, u_{12}, u_{21}, T_{12}$ and $T_{21}$ will be determined using conservation of \textcolor{black}{the number of particles,} total momentum and \textcolor{black}{energy} \eqref{cons2}, together with some symmetry considerations.
If we assume that \begin{align} n_{12}=n_1 \quad \text{and} \quad n_{21}=n_2,  
\label{density} 
\end{align}
we have conservation of the number of particles, see Theorem 2.1 in \cite{Pirner}.
If we further assume that $u_{12}$ is a linear combination of $u_1$ and $u_2$
 \begin{align}
u_{12}= \delta u_1 + (1- \delta) u_2, \quad \delta \in \mathbb{R},
\label{convexvel}
\end{align} then we have conservation of total momentum
provided that
\begin{align}
u_{21}=u_2 - \frac{m_1}{m_2} \varepsilon (1- \delta ) (u_2 - u_1),
\label{veloc}
\end{align}
see Theorem 2.2 in \cite{Pirner}.
If we additionally assume that $T_{12}$ is of the following form
\begin{align}
\begin{split}
T_{12} &=  \alpha T_1 + ( 1 - \alpha) T_2 + \gamma |u_1 - u_2 | ^2,  \quad 0 \leq \alpha \leq 1, \gamma \geq 0 ,
\label{contemp}
\end{split}
\end{align}
then we have conservation of total energy
provided that
\begin{align}
\begin{split}
T_{21} =\left[ \frac{1}{3} \varepsilon m_1 (1- \delta) \left( \frac{m_1}{m_2} \varepsilon ( \delta - 1) + \delta +1 \right) - \varepsilon \gamma \right] |u_1 - u_2|^2 \\+ \varepsilon ( 1 - \alpha ) T_1 + ( 1- \varepsilon ( 1 - \alpha)) T_2,
\label{temp}
\end{split}
\end{align}
see Theorem 2.3 in \cite{Pirner}.
In order to ensure the positivity of all temperatures, we need to restrict $\delta$ and $\gamma$ to 
 \begin{align}
0 \leq \gamma  \leq \frac{m_1}{3} (1-\delta) \left[(1 + \frac{m_1}{m_2} \varepsilon ) \delta + 1 - \frac{m_1}{m_2} \varepsilon \right],
 \label{gamma}
 \end{align}
and
\begin{align}
 \frac{ \frac{m_1}{m_2}\varepsilon - 1}{1+\frac{m_1}{m_2}\varepsilon} \leq  \delta \leq 1,
\label{gammapos}
\end{align}
see Theorem 2.5 in \cite{Pirner}. For all these choices one can prove the entropy inequalities \eqref{entropy}, see Theorem 2.7 in \cite{Pirner}. We observe that we have free parameters $\alpha, \delta, \gamma$.
We keep the free parameters to be as general as possible. We will discuss the meaning and possible choices in section \ref{sec:parameters}.

\section{Theoretical results of this model}
\label{sec3}
In this section, we give an overview over recent theoretical results for the model presented in section \ref{sec2} concerning existence of solutions and large-time behaviour. To start with, one can prove an existence and uniqueness result of mild solutions in the periodic setting in space under certain conditions on the initial data and the collision frequencies. The proof is presented in \cite{Pirner3}. Another existence result concerning the existence of a unique global-in-time classical solution when the initial data perturbed slightly from a global equilibrium can be found in \cite{Koreaner}.

Moreover, one can prove the following results on the large-time behaviour \cite{Crestetto_2012}.
 We denote the entropy of a function $f$  by $H(f)= \int f \ln f dv$ and the relative entropy of $f$ and $g$ by $H(f|g)= \int f \ln \frac{f}{g} dv$.
Then one can prove the following results on the large-time behavior  \cite{Crestetto_2012}.
\begin{theorem}
 Suppose  that $\nu_{12}$ is constant in time. Then, in the space homogeneous case 
we have the following decay rate of the distribution functions $f_1$ and $f_2$
$$ || f_k - M_k ||_{L^1(dv)} \leq 4 e^{- \frac{1}{2} Ct } [ H(f_1^0|M_1^0) + H(f_2^0 | M_2^0)]^{\frac{1}{2}}, \quad k=1,2,$$
where $C$ is the constant given by 
$$ C= \min \lbrace \nu_{11} n_1 + \nu_{12} n_2, \nu_{22} n_2 + \nu_{21} n_1 \rbrace.$$
and the index 0 denotes the value at time t = 0.
\end{theorem}
  \begin{theorem}\label{th:estimate_vel}
Suppose that $\nu_{12}$ is constant in time. In the space-homogeneous case, we have the following relaxation rate
\begin{align}
    \partial_t(u_1-u_2) = \nu_{12} (1- \delta) (n_2+ \frac{m_1}{m_2} n_1) (u_2-u_1)
    \label{rel_vel}
\end{align}
and a decay rate of the mean velocities
\begin{equation*}
|u_1(t) - u_2(t)|^2 = e^{- 2 \nu_{12} (1- \delta)\left(n_2+\frac{m_1}{m_2} n_1\right) t} |u_1(0) - u_2(0)|^2.
\end{equation*}
\end{theorem}
  \begin{theorem}\label{th:estimate_temp} 
 Suppose $\nu_{12}$ is constant in time. In the space-homogeneous case, we have the following relaxation rate
 \begin{align}
     \partial_t(T_1-T_2) = - C_1 (T_1-T_2) + C_2 |u_1-u_2|^2
     \label{rel_temp}
 \end{align}
and a decay rate of the temperatures
\begin{equation*}
\begin{split}
\textcolor{black}{T_1(t) - T_2(t) = 
 e^{- C_1t} \left[T_1(0) - T_2(0)+\frac{C_2}{C_1-C_3 } ( e^{(C_1-C_3) t} - 1) |u_1(0) - u_2(0)|^2 \right],}
 \end{split}
\end{equation*}
where the constants are defined by
\begin{align*}
C_1&=(1- \alpha)\nu_{12}\left(n_2+n_1\right),\\
C_2&=\nu_{12}\left(n_2\left((1-\delta)^2+\frac{\gamma}{m_1}\right)-n_1\left(1-\delta^2-\frac{\gamma}{m_1}\right)\right),\\
C_3&=2 \nu_{12} (1- \delta)\left(n_2+\frac{m_1}{m_2} n_1\right).
\end{align*}
\end{theorem}

There are also results in the space-inhomogeneous case for the linearized collision operator, see \cite{Pirner_Liu}. In their article, the authors study hypocoercivity for the linearized BGK model for gas mixtures in continuous phase space. By constructing an entropy functional, they prove exponential relaxation to global equilibrium with explicit rates. The strategy is based on the entropy and spectral methods adapting Lyapunov's direct method as presented in \cite{achleitnerlinear} for the one species linearized BGK model. 
\section{Possible choices and meaning of the free parameters} \label{sec:parameters}
In this section, we deal with the meaning and possible choices of the free parameters. One possibility is that we can choose the parameters such that we can generate special cases in the literature \cite{gross_krook1956,hamel1965,asinari,Garzo1989,Sofonea2001,Cercignani_1975,Greene,Bobylev,haack}. For instance if we choose  $\varepsilon=1$, $ \delta= \frac{m_1}{m_1+m_2}$, $\alpha=\frac{m_1^2+m_2^2}{(m_1+m_2)^2}$ and $\gamma=\frac{m_1 m_2}{(m_1+m_2)^2} \frac{m_2}{3}$, we obtain the model by Hamel in \cite{hamel1965}. 

Another possibility  is to choose the parameters in a way such that the macroscopic exchange terms of momentum and energy can be matched in a certain way for example that they coincide with the ones for the Boltzmann equation. For this, we first present the macroscopic equations with exchange terms of the BGK model \eqref{BGK}.
If we multiply the BGK model for gas mixtures by $1, m_j v, m_j \frac{|v|^2}{2}$ and integrate with respect to $v$, we obtain the following macroscopic conservation laws
\begin{multline*}
\\
\partial_t n_1 + \nabla_x \cdot (n_1 u_1)=0, 
\\
\partial_t n_2 + \nabla_x \cdot (n_2 u_2)=0,
\\
 \partial_t(m_1 n_1 u_1)+\nabla_x \cdot \int m_1 v \otimes v f_1(v) dv + \nabla_x \cdot (m_1 n_1 u_1 \otimes u_1  )  =  f_{m_{1,2}},
\\
 \partial_t(m_2 n_2 u_2)+\nabla_x \cdot \mathbb{P}_2 + \nabla_x \cdot (m_2 n_2 u_2 \otimes u_2  )  = 
f_{m_{2,1}},
\\
\partial_t \left(\frac{m_1}{2} n_1 |u_1|^2 + \frac{3}{2} n_1 T_1 \right) + \nabla_x \cdot \int m_1 |v|^2 v f(v) dv  =  F_{E_{1,2}},
\\
\partial_t \left(\frac{m_2}{2} n_2 |u_2|^2 + \frac{3}{2} n_2 T_2 \right) + \nabla_x \cdot Q_2  = F_{E_{2,1}},
\\
\end{multline*}
with exchange terms $f_{m_{i,j}}$ and $F_{E_{i,j}}$ given by
\begin{align*}
f_{m_{1,2}}&= - f_{m_{2,1}} = m_1 \nu_{12} n_1 n_2 (1 - \delta) (u_2 - u_1), \\
F_{m_{1,2}}&= - F_{m_{2,1}} \\& = \left[\nu_{12} \frac{1}{2} n_1 n_2 m_1 (\delta -1) (u_1 + u_2 + \delta(u_1-u_2)) + \frac{1}{2} \nu_{12} n_1 n_2 \gamma (u_1 - u_2) \right] \cdot (u_1-u_2)\\ &+ \frac{3}{2} \varepsilon \nu_{21} n_1 n_2 (1-\alpha) (T_2-T_1).\label{macrosequ1}
\end{align*}

Here, we can observe a physical meaning of $\alpha$ and $\delta$. We see that $\alpha$ and $\delta$ show up in the exchange terms of momentum and energy as parameters in front of the relaxation of $u_1$ towards $u_2$ and $T_1$ towards $T_2$. So they determine, together with the collision frequencies, the speed of relaxation of the mean velocities and the temperatures to a common value. This can already be observed in Theorem \ref{th:estimate_vel} and Theorem \ref{th:estimate_temp}.

Next we follow Chapter 4.1 in \cite{haack} and compare the relaxation rates in the space-homogeneous case to the relaxation rates for the space-homogeneous Boltzmann equation. In \cite{haack}, they find values for $\nu_{kj}$ such that either the relaxation rate for the mean velocities \eqref{rel_vel} or the relaxation for the temperatures \eqref{rel_temp} coincides with the corresponding rate of the Boltzmann equation. But using the free parameters $\alpha$, $\delta$ and $\gamma$ we are able to match both of the relaxation rates at the same time. For this, we compare the coefficients of the terms $u_2 - u_1$, $T_2 - T_1$ and $\vert u_2 - u_1 \vert^2$ in these Boltzmann relaxation rates and the BGK relaxation rates \eqref{rel_vel} and \eqref{rel_temp}, and we derive the values of the parameters for this model:
\begin{align*}
(u_2 - u_1) \text{-term:} \quad \delta &= 1- \frac{\alpha_{12}}{\nu_{12}} \frac{m_1 + m_2}{2  } \frac{m_1 n_1 + m_2 n_2}{m_1 n_1 m_2 n_2} \left(  n_1 \frac{m_1}{m_2} + n_2 \right)^{-1},\\
(T_2 - T_1) \text{-term:} \quad \alpha &= 1- \frac{\alpha_{12}}{\nu_{12} n_2 n_1} , \\
\vert u_2 - u_1 \vert^2 \text{-term:} \quad \gamma &= \frac{1}{3}\left(  n_1  + n_2 \right)^{-1} \left[\frac{\alpha_{12}}{\nu_{12}} \frac{m_2 n_2 - m_1 n_1}{ n_2 n_1 } - m_1 n_2 (1-\delta)^2 +  m_1 n_1 (1-\delta^2)   \right],
\end{align*} 
where $\alpha_{12}$ is a coefficient for energy transfer coming from Boltzmann equation, see \cite{haack} and references therein. Additionally, the constraints \eqref{contemp}, \eqref{gamma} and \eqref{gammapos} need to be satisfied. This can be verified by a corresponding choice of $\nu_{kj}$. One possibility is
\begin{align}
    \nu_{kj} = \frac{1}{2} \frac{\alpha_{kj}}{n_k n_j} \frac{(m_k+m_j)^2}{m_k m_j}
\end{align}
and for $1\geq \varepsilon=\frac{m_j}{m_k}$ (cf. in a plasma).

\section{On existing numerical schemes}
\label{sec4}
In the literature, \textcolor{black}{various approaches for the discretization of kinetic equations can be found}, including schemes for the one-species BGK model. \textcolor{black}{Contributions} in numerics for multi-species BGK models have strongly increased in the last years. 

To start with, we give a short overview over existing numerical methods for the one-species BGK equation. Since the contributions are very crowded, we do not claim completeness. Many ideas can be carried over to the \textcolor{black}{discretization of} multi-species BGK equations, \textcolor{black}{and we} conclude with identified publications on numerical schemes for multi-species BGK equations which can be written in the form \eqref{BGK}.\\
\\
The (one-species) Boltzmann equation \textcolor{black}{captures physical phenomena} very well at the kinetic level \cite{Cercignani_1975}. Nevertheless, numerical computation is expensive. The fastest algorithms for evaluating the Boltzmann collision operator are spectral methods with special kernels \cite{Mouhot2006}. This motivates the BGK equation as an approximation of the Boltzmann equation: Even though the dimensionality is \textcolor{black}{as high as for the Boltzmann equation}, the \textcolor{black}{BGK} interaction term is better to handle and explicitly computable \cite{Puppo_2007,Jin_2010}. Hence, the computational cost is much less compared to the Boltzmann equation while maintaining most of the physical properties. \cite{Struchtrupp, MieussensStruchtrup, Holway}

The computational advantages are also useful for penalization techniques \cite{Jin_2010} where the BGK equation is solved as preconditioner for the numerical solution of the Boltzmann equation. This idea is generalized to the multi-species setting in \cite{JinLi2013}. In \cite{DegondDimarcoPareschi2011}, \textcolor{black}{the authors develop an improved Monte Carlo method for the BGK equation. This is supposed to be a first step towards} an improved Monte Carlo simulation of the Boltzmann equation. Moreover, the BGK \textcolor{black}{approach is useful when coupling different domains in which} the regimes range from equilibrium to very rarefied \cite{AlaiaPuppo2012}. \\
\\
\textcolor{black}{A fully-discrete scheme requires the discretization in (microscopic) velocity, space and time.}
First we consider the discretization in velocity before we look at the space and time variables.

Having the microscopic velocities as independent variables introduces both more degrees of freedom and more difficulties. Due to the high dimensionality, it is recommendable to use coarse grids \cite{Puppo2018} which then poses challenges regarding errors in the macroscopic quantities. This can be tackled when the conservation properties \eqref{cons1} are fulfilled at the discrete level. The handling of discrete moments, a discrete entropy and the corresponding discrete Maxwellians is discussed in \cite{Mieussens2000}. Another approach to \textcolor{black}{fulfil the conservation laws at the discrete level} is given by a constrained $L^2$-projection in \cite{Gamba2009}. 

Being interested in macroscopic quantities only, the Chu reduction is a possible approach to lower the dimensionality if there are more degrees of freedom in velocity than in space \cite{Chu}. \textcolor{black}{Using the Chu reduction, one follows} the evolution of appropriate integrals of the distribution functions, but these integrals \textcolor{black}{do not correspond to} macroscopic quantities yet. This method reduces the computational costs considerably. 

When the mean velocities $u(x,t)$ cover a wide range and small temperatures are encountered, grid adaption becomes an important tool. This issue is tackled more and more in the last decade, e.g. in \cite{BrullMieussens2014,BernardIolloPuppo2014,BoscarinoChoRussoPreprint2021,HittingerBanks2013}.  

\textcolor{black}{This leads us to} another advantage of multi-species BGK equations. \textcolor{black}{In case of the multi-species Boltzmann equations, a large mass ratio of the species, which results in very distinct thermal speeds, requires an expensive grid resolution \cite{MunafoTorresHaackGambaMagin2014}. As particles of different species only interact through moments in the BGK model, the evolution of each species can be numerically solved on separate grids \cite{haack} which might be an important ingredient for an efficient simulation. } \\
\\
There can be found many different approaches for the space discretization in the literature as the BGK equation shares the same transport term with many other kinetic equations such as Boltzmann, Fokker-Planck, Vlasov, etc.

The transport term being hyperbolic, a finite volume discretization is often used. High-resolutions can be obtained by weighted essentially non-oscillatory (WENO) or discontinuous Galerkin (DG) schemes.  However, for orders higher than two \textcolor{black}{the corresponding formulation of the relaxation term requires additional care} because it does not suffice to consider the relaxation of the cell averages, but the cell averages of the relaxation term need to be calculated. \cite{Mieussens2000,Puppo_2007,HuJinLi2017,AyusoCarrilloShu2011,ChengGambaProft2011}

Another \textcolor{black}{convenient choice is the semi-Lagrange method. The characteristics are followed exactly} which requires an interpolation for the evaluation of the corresponding foot point. By conservative reconstructions or corrections, these methods can be kept conservative also for higher orders \cite{SonnendrueckerRocheBertrandGhizzo1998,CrouseillesMehrenbergerSonnendruecker2010,Dimarco_2014,ChoBoscarinoRussoYun2021,QiuShu2011}.

In \cite{DimarcoLoubere2012}, the authors present \textcolor{black}{an} approach for an efficient scheme based on discrete velocity models and semi-Lagrangian methods. \textcolor{black}{In contrast to standard semi-Lagrangian schemes}, the distribution function needs not to be reconstructed at each time step which of course accelerates numerical computations.

For the Vlasov equation, the most used method is the Particle In Cell (PIC) method \cite{FilbetSonnendrueckerBertrand2001}. But to our knowledge, it is less used for equations with interaction terms when hydrodynamic effects become more important. In \cite{FilbetSonnendrueckerBertrand2001}, the authors shortly discuss different methods for the Vlasov equation and then introduce their positive and flux conservative method (PFC). \\
\\
For the interaction term, a time implicit formulation is \textcolor{black}{often} chosen since the right-hand side becomes stiff when the collision frequencies become large (close to the hydrodynamic regime). By the implicit discretization, one can avoid tiny time steps coming from stability issues. Thanks to the special structure of the interaction term, the implication is comparably easy manageable, and the equation stays explicitly solvable \cite{Puppo_2007,Jin_2010}.  

Usually, the transport part is evaluated explicitly. It is combined with the interaction term e.g. by splitting methods \cite{CoronPerthame1991,haack2}, implicit-explicit Runge-Kutta (IMEX RK) schemes \cite{PareschiRusso2005,Puppo_2007} or IMEX multistep methods \cite{DimarcoPareschi2017}. 

Splitting methods must be treated with care when the right-hand side becomes stiff. In \cite{Jin1995}, the author shows that the (second-order) Strang splitting reduces actually to a first-order approximation of the equilibrium equation in the hydrodynamic limit. This leads us to so-called asymptotic-preserving (AP) schemes which provide an adequate discretization also of the limiting equations. \textcolor{black}{Using AP schemes}, the correct equilibrium solutions are preserved \cite{Puppo2018,HuJinLi2017,Bernard_2015,Dimarco_2014,Jin_2010}. This issue is addressed more and more in the recent years. In this context, we also want to mention the micro-macro decomposition and the parity decomposition/AP splitting. For the \textcolor{black}{former} approach, the distribution function is written as a sum of its equilibrium (macro) part and the \textcolor{black}{remnant} which represents the kinetic (micro) part. This results in one microscopic and one macroscopic equation which can be solved by individual and adequate methods \cite{Crestetto_2013}. For the latter approach, the distribution function is decomposed by an even and an odd parity. A new system of equations can be derived with only one time scale where splitting techniques can be applied  \cite{JinPareschi2000,Dimarco_2014}. \\
\\
More physics is (re)introduced by multi-species models. At the discrete level, many ideas can be carried over from the single-species schemes. In the following, we give contributions of numerical schemes for multi-species BGK equations, which can be written in the form of \eqref{BGK}.

In \cite{Crestetto_2012}, \textcolor{black}{the authors} extend the work of \cite{Crestetto_2013} for the multi-species model. They perform a micro-macro decomposition and focus on the fluid limit. The micro part is solved by a particle method, whereas the macro part (depicting the fluid part) is solved by a standard finite volume approach. Here, an additional force term with an electric field is considered.

In \cite{haack}, \textcolor{black}{the authors} are interested in capturing physical transport coefficient. They use the additional degrees of freedom in the multi-species setting to match relaxation rates in the space homogeneous case equivalent to the Boltzmann ones.  An extension to space inhomogeneous simulations is done in \cite{haack2} where they additionally examine the coupling to electric fields.

In \cite{BoscarinoChoGropppiRusso}, \textcolor{black}{the authors} compare numerical results for different multi-species BGK models, where one of these models is a special case of \eqref{BGK}.

\textcolor{black}{A BGK model for gas mixtures is extended to velocity dependent collision frequencies in \cite{Pirner_velocity}. Collision frequencies influence the relaxation process and the resulting hydrodynamic behavior such that they also become important when calculating transport coefficients. The class of models in \cite{Pirner_velocity} captures a model of the form \eqref{BGK} as a special case, but in general the Maxwellians are replaced by more sophisticated Gaussian functions.
Numerical schemes for this kind of equations have been developed in \cite{Sandra_Numerik}. The key new ingredient is a solver based on a convex entropy minimization problem which makes possible an implicit treatment of the BGK operator. }


\section{Conclusions and Outlook}
We recalled a general BGK model for gas mixtures with two collision terms for two species and its theoretical properties. \textcolor{black}{We presented results on the H-Theorem, existence of solutions and large-time behaviour. }The special feature of the presented model are the free parameters $\alpha$, $\delta$ and $\gamma$. They influence the exchange of momentum and energy and can be set such that the model's behavior matches the physics or coincides with another (more specialized) model.

An overview over existing literature on numerics for this kind of models was given.

\textcolor{black}{ However, BGK- type models often lack on correct parameters in the continuum limit like the Prandtl number. Therefore these models can be used as a basis for more extended models like ES-BGK models or BGK models with velocity dependent collision frequency. As a future work the Chapman-Enskog expansion of 
such
models can be computed. Then the transport coefficients of all these models can be compared and eventually extended to match all parameters in the macroscopic equations. Here, the free parameters in the BGK model for monoatomic molecules with a sum of interaction terms might be useful.  }


%
%

\section*{Acknowledgement}
Marlies Pirner was supported by the Alexander von Humboldt foundation.

\end{document}